\documentclass[12pt,reqno]{amsart}

\usepackage{amsmath,amsthm}
\usepackage{amssymb}
\usepackage{euscript}

\numberwithin{equation}{subsection}

\newtheorem{Thm}[subsection]{Theorem}
\newtheorem{Lem}[subsection]{Lemma}
\newtheorem{Prop}[subsection]{Proposition}
\newtheorem{Cor}[subsection]{Corollary}

\theoremstyle{definition}

\newtheorem{Exm}[subsection]{Example}
\newtheorem{Rem}[subsection]{Remark}

\begin{document}
\title[]{On trivialities of Chern classes}
\author{Aniruddha C. Naolekar}
\address{Indian Statistical Institute, 8th Mile, Mysore Road, RVCE Post, Bangalore 560059, INDIA.}
\email{ani@isibang.ac.in}

\author{ Ajay Singh Thakur}
\address{Department of Mathematics, University of Haifa, Mount Carmel, Haifa 31905, ISRAEL.}
\email{thakur@math.haifa.ac.il}

\begin{abstract} A finite $CW$-complex $X$ is $C$-trivial if for every complex vector bundle $\xi$ over $X$, the total Chern class $c(\xi)=1$. In this note we completely determine when each of the following spaces are $C$-trivial: suspensions of stunted real projective spaces, suspensions of stunted complex projective spaces and suspensions of stunted quaternionic projective spaces. 
\end{abstract}


\email{}

\date{}
\maketitle

\section{Introduction}

A $CW$-complex $X$ is said to be $C$-trivial if for any complex vector bundle $\xi$ over $X$, the total Chern class $c(\xi)=1$. 

A related notion is that of $W$-triviality. A $CW$-complex $X$ is said to be $W$-trivial if for any real vector bundle $\eta$ over $X$, the total Stiefel-Whitney class $w(\eta)=1$. A central result in this direction is a theorem of Atiyah-Hirzebruch (\cite{atiyahirz}, Theorem\,2) which states that for any finite $CW$-complex $X$, the $9$-fold suspension $\Sigma^9X$ of $X$ is always $W$-trivial. In particular, the spheres $S^k$ are all $W$-trivial for $k\geq 9$. In fact, a sphere $S^k$ is $W$-trivial if and only if $k\neq 1,2,4,8$ (see \cite{atiyahirz}, Theorem\,1).

Understanding which spaces 
are $W$-trivial has been of some interest in recent times (see \cite{aniajay}, \cite{tanaka}, \cite{ajay} and the references therein). 
In \cite{tanaka}, the author has completely determined which suspensions $\Sigma^k \mathbb F\mathbb P^n$ are $W$-trivial. Here $\mathbb F$ denotes either the field $\mathbb R$ of real number, the field $\mathbb C$ of complex numbers or the skew-field $\mathbb H$ of quaternions and $\mathbb F\mathbb P^n$ denotes the appropriate projective space. In \cite{ajay}, the second named author has determined, in most cases, which suspensions of the Dold manifolds are $W$-trivial. In \cite{aniajay}, the authors have completely determined which suspensions of the stunted real projective spaces are $W$-trivial. 

In this note we study the notion of $C$-triviality and determine whether some of the familiar spaces and their suspensions are $C$-trivial. 
To begin with it is well known that there is no analogue of the Atiyah-Hirzebruch theorem for Chern classes. Indeed, by the Bott integrality theorem (see Theorem\,\ref{bott} below for the precise statement) the even dimensional spheres are not $C$-trivial. 
Thus if $d>4$, then $S^{2d}$ is $W$-trivial but not $C$-trivial. The circle $S^1$ is $C$-trivial but not $W$-trivial. 
We give other examples in the sequel. 
However, there are sufficient conditions under which one implies the other. We point out conditions under which this happens. 

In this note we completely determine when the suspension of a stunted real projective space is $C$-trivial (see Theorem\,\ref{secondtheorem} and Theorem\,\ref{thirdtheorem} below). We also 
completely determine which suspensions of the stunted complex and quaternionic projective spaces are $C$-trivial (see Corollary\,\ref{complex} below). 
 
We now state the main results of this paper.  
The following theorem completely describes which suspensions $\Sigma^k\mathbb R\mathbb P^n$ of the real projective spaces are $C$-trivial. Since $S^1 =  \mathbb R \mathbb P^1$ is $C$-trivial and $\mathbb R \mathbb P^n$ is not $C$-trivial for $n >1$,  we shall assume $k>0$.

\begin{Thm} \label{secondtheorem} Let $X^k_n=\Sigma^k \mathbb R\mathbb P^n$ with $k,n>0$.  Then $X^k_n$ is not $C$-trivial if and only if one of the following conditions is satisfied. 
\begin{enumerate}
\item $k,n$ are both odd. 
\item $k=2,4$ and $n\geq k$.
\end{enumerate}
\end{Thm}


%
%
Next we look at the suspensions of the stunted real projective spaces.  To state the results we introduce the following notations. 
Let $X_{m,n}$ denote the stunted real projective space 
$$X_{m,n}=\mathbb R\mathbb P^m/\mathbb R\mathbb P^n$$
and $X^k_{m,n}$ the $k$-fold suspension 
$$X^k_{m,n}=\Sigma^k\left(\mathbb R\mathbb P^m/\mathbb R\mathbb P^n\right).$$

\begin{Thm}\label{thirdtheorem} 
Let $X^k_{m,n}$ be as above with $k\geq 0$ and $0<n<m$. 
\begin{enumerate}
\item If $k$ is odd and $m$ is even, then $X^k_{m,n}$ is $C$-trivial. 
\item If $k,m$ are odd, then $X^k_{m,n}$ is not $C$-trivial.
\item If $ n= 2t$, then $X_{m,n}$ is $C$-trivial if and only if $m< 2^{t+1}$. 
\item If $k,n$ are even and $k\geq 2$, then $X^k_{m,n}$ is $C$-trivial. 
\item If $k$ is even and $n$ is odd, then $X^k_{m,n}$ is not $C$-trivial. 
\end{enumerate}
\end{Thm}

The paper is organized as follows. In Section 2 we prove some general facts about $C$-triviality. In Section 3 we determine which suspensions of stunted complex and quaternionic projective spaces are $C$-trivial and prove the main theorems. 

{\em Conventions.} By a space we mean a finite connected $CW$-complex. Given a map $\alpha:X\longrightarrow Y$ between spaces the induced homomorphism in $K$-theory and singular cohomology will again be denoted by $\alpha$. 


\section{Generalities}

In this section we prove some general facts about $C$-triviality. We give sufficient conditions under which $W$-triviality implies $C$-triviality and 
$C$-triviality implies $W$-triviality. 

To begin, note that a finite $CW$-complex $X$ is $C$-trivial if the reduced complex $\widetilde{K}$-group $\widetilde{K}(X)=0$.  

\begin{Lem} \label{prelim} For any space $X$, we have the following. 
\begin{enumerate}
\item If $H^{2s}(X;\mathbb Z)=0$ for all $s> 0$, then $X$ is $C$-trivial. 
\item If $X$ has cells only in odd dimensions, then $X$ is $C$-trivial. 
\item If $X$ is $C$-trivial, then $H^2(X;\mathbb Z)=0$
\item If $X$ is $W$-trivial and the mod-$2$ reduction homomorphism $\rho_2:H^{2i}(X;\mathbb Z)\longrightarrow H^{2i}(X;\mathbb Z_2)$ is monomorphic  for all $i>0$, then $X$ is $C$-trivial. 
\end{enumerate}
\end{Lem}
{\bf Proof.} We omit the easy proofs of (1)-(3). The claim  (4) follows from the fact that for a complex bundle $\xi$ we have $\rho_2(c_i(\xi))=w_{2i}(\xi_{\mathbb R})$. Here $\xi_{\mathbb R}$ denotes the underlying real bundle of $\xi$. \qed

Thus the real and complex projective spaces $\mathbb R\mathbb P^n$ ($n>1$) and $\mathbb C\mathbb P^n$ are not $C$-trivial as their second integral cohomology group is non-zero. 

The following observations are straightforward and we omit their easy proofs. 

\begin{Lem}\label{prelim1} 
Let $f:X\longrightarrow Y$ be a map between spaces. 
\begin{enumerate}
\item If $f:\widetilde{K}(Y)\longrightarrow \widetilde{K}(X)$ is onto and $Y$ is $C$-trivial, then $X$ is $C$-trivial. 
\item Suppose $f:H^{2i}(Y;\mathbb Z)\longrightarrow H^{2i}(X;\mathbb Z)$ is a monomorphism for all $i>0$. Then if $X$ is $C$-trivial so is $Y$. \qed
\end{enumerate}
\end{Lem}

We shall use the above observations in the sequel sometimes without an explicit reference. We have already noted that the even dimensional spheres  $S^{2d}$ with $d>4$ are examples of spaces that are $W$-trivial but not $C$-trivial and also that $S^1$ is $C$-trivial but not $W$-trivial. 
We give some more examples below.

\begin{Exm}
Let $X=\Sigma\mathbb R\mathbb P^2$ be the suspension of the real projective space $\mathbb R\mathbb P^2$. Then $X$ has non-zero integral cohomology  only in degree $3$ and hence $X$ must be $C$-trivial. It is known that 
$X$ is not $W$-trivial (\cite{tanaka}, Theorem 1.4).    
\end{Exm} 


\begin{Exm}
Let $X=M(\mathbb Z_3,1)$ be the Moore space of type $(\mathbb Z_3,1)$. This is a $2$-dimensional $CW$-complex. Since the second (integral) cohomology of $X$ is non-zero, $X$ is not $C$-trivial. However, as $H^i(X;\mathbb Z_2)=0$ 
for $i>0$, $X$ is $W$-trivial. 
\end{Exm}




We now state a necessary condition for a space to be $C$-trivial. 

\begin{Lem}
Let $X$ be a $C$-trivial space. Then for any real bundle $\xi$ over $X$ we have $w_i^2(\xi)=0$ for all $i>0$.  
\end{Lem}
{\bf Proof.} Let $\eta$ denote the underlying real bundle of the complexification $\xi\otimes \mathbb C$ of $\xi$. Then the total Stiefel-Whitney class 
of $\eta$ is given by 
$$w(\eta)=1+w_1^2+w_2^2+\cdots$$
where $w_i=w_i(\xi)$. If $w(\eta)\neq 1$, then clearly $c(\xi\otimes \mathbb C)\neq 1$. This completes the proof. \qed

In particular, if there exists a real bundle $\xi$ over $X$ with $w_i^2(\xi)\neq 0$ for some $i>0$, then $X$ cannot be $C$-trivial. 
The above lemma, in particular, implies that the quaternionic projective space $\mathbb H\mathbb P^n$ is not $C$-trivial for $n>1$. This is because for the canonical line bundle $\xi$ over $\mathbb H\mathbb P^n$ we have $w_4^2(\xi)\neq 0$. That $\mathbb H\mathbb P^1=S^4$ is  not $C$-trivial is clear. The (even dimensional) spheres show that the converse of the above lemma is not true. 

For a space $X$, we have the realification homomorphism $r:\widetilde{K}(X)\longrightarrow \widetilde{KO}(X)$. 
The following lemma gives a sufficient condition for a $C$-trivial space to be $W$-trivial. 
 
\begin{Lem}
Suppose that $r:\widetilde{K}(X)\longrightarrow \widetilde{KO}(X)$ is onto. If $X$ is $C$-trivial, then $X$ is $W$-trivial. 
\end{Lem}
{\bf Proof.} The surjectivity of $r$ implies that every real vector bundle $\xi$ over $X$ is stably equivalent to the underlying real bundle 
$\eta_{\mathbb R}$ of a complex vector bundle $\eta$ over $X$. Now $w(\xi) =w(\eta_{\mathbb R})=1$. This completes the proof. \qed

We mention that the converse to the above lemma is not true. Indeed, if $X=S^6$, then as $\widetilde{KO}(S^6)=0$, $X$ is $W$-trivial but $X$ is not 
$C$-trivial. 

Some of our proofs depend upon the following important observation which gives 
a sufficient condition for a $W$-trivial space to be $C$-trivial. 

\begin{Prop}\label{onetytwo}
Let $X$ be a space. Assume that $H_*(X;\mathbb Z)$ is concentrated in odd degrees and is direct sum of copies of $\mathbb Z _2$. If $X$ is $W$-trivial, then $X$ is $C$-trivial. 
\end{Prop}
{\bf Proof.} Given the assumptions on the integral homology of $X$ it follows from the universal coefficients theorem that the integral cohomology of $X$ is concentrated in even degrees and is a direct sum of copies of $\mathbb Z_2$. Again, by the universal coefficients theorem, it follows that $H^{2i}(X;\mathbb Z_2)$ is a direct sum of copies of $\mathbb Z_2$ with the same number of $\mathbb Z_2$ summands as in the integral cohomology. 
As $X$ has no integral cohomology in odd degrees, the mod-$2$ reduction map 
$$\rho_2:H^{2i}(X;\mathbb Z)\longrightarrow H^{2i}(X;\mathbb Z_2)$$
is surjective and hence an isomorphism as both the groups are a direct sum of equal number of copies of $\mathbb Z_2$. The proposition now follows from Lemma\,\ref{prelim} (4). This completes the proof. \qed

\begin{Rem}
We remark that the converse of the above proposition is not true. For consider the space $X^6_2=\Sigma^6\mathbb R\mathbb P^2$. Then $X^6_2$ satisfies the conditions of the above proposition. By Theorem\,\ref{secondtheorem}, $X^6_2$ is $C$-trivial. However $X^6_2$ is not $W$-trivial, by Theorem\,1.4 (3) of \cite{tanaka}. 
\end{Rem}

We end this section by noting that second suspension of a $C$-trivial space is $C$-trivial. The idea of the proof is  same 
as that of Theorem\,1 in \cite{atiyahirz} and Theorem\,1.1 of \cite{tanaka}. 

\begin{Thm}\label{hahaha}
Let $X$ be a $C$-trivial space. Then the second suspension $\Sigma^2X$ of $X$ is $C$-trivial. 
\end{Thm}
{\bf Proof.} Let $\pi_1,\pi_2$ denote the projection maps of $S^2\times X$ to the first and second factor respectively and $p:S^2\times X\longrightarrow \Sigma^2X$ the quotient map. Let $\nu$ denote the Hopf bundle over $S^2$. Given a bundle $\xi$ over $\Sigma^2X$ there exists, by Bott periodicity, a bundle $\theta$ of rank $n$ (say) over $X$ such that $p^*\xi$ is stably isomorphic to the tensor product 
$$(\pi_1^*\nu-1)\otimes (\pi_2^*\theta-n).$$
Thus 
$$c(p^*\xi)=c(\pi_1^*\nu\otimes\pi_2^*\theta)c(\pi_1^*\nu)^{-n}c(\pi_2^*\theta)^{-2}.$$
Since $X$ is $C$-trivial we have $c(\pi_2^*\theta)=1$. It is clear that $c(\pi_1^*\nu)=1+t\times 1$ where $t\in H^2(S^2;\mathbb Z)$ is a generator. 
Finally, one checks that 
$$c(\pi_1^*\nu\otimes\pi_2^*\theta)=(1+t\times 1)^n$$
so that 
$$c(p^*\xi)=1.$$
But as $f:H^*(\Sigma^2X;\mathbb Z)\longrightarrow H^*(S^2\times X;\mathbb Z)$ is a monomorphism it follows that $c(\xi)=1$. This completes the proof. \qed

\section{Proof of Theorems\,\ref{secondtheorem} and \ref{thirdtheorem}}

In this section we prove the main theorems and derive some consequences. We first state the Bott integrality theorem which we shall use in the sequel. 

\begin{Thm} {\rm(Bott integrality theorem)}\label{bott} {\rm(\cite{hus}, Chapter 20, Corollary\,$9.8$)} Let $a\in H^{2n}(S^{2n};\mathbb Z)$ be a generator. For any complex vector bundle $\xi$ over $S^{2n}$, the Chern class $c_n(\xi)$ is divisible by $(n-1)!a$. For each $m$ divisible by $(n-1)!$, there exists a unique $\xi\in\widetilde{K}(S^{2n})$ with $c_n(\xi)=ma$. \qed
\end{Thm}

The Bott integrality theorem implies that if $\xi$ is a complex vector bundle over the even dimensional sphere $S^{2n}$ with $n \geq 3$, then $c_n(\xi) \in H^{2n}(S^{2n};\mathbb Z)$ is an even multiple of a generator.

The following observation is now immediate from the Bott integrality theorem. 

\begin{Prop}\label{firstcor}
Let $X$ be a finite $CW$-complex. Assume that there is a map $\alpha:X\longrightarrow S^d$ such that the homomorphism 
$\alpha:H^d(S^d;\mathbb Z)\longrightarrow H^d(X;\mathbb Z)$ is injective. Then $\Sigma^kX$ is not $C$-trivial whenever $k+d$ even. 
\end{Prop}
{\bf Proof.} The map $\Sigma^kf$ induces a monomorphism in cohomology in degree $(k+d)$ for all $k$. If $(k+d)$ is even, then as there exists a vector bundle $\xi$ over $S^{k+d}$ with $c(\xi)\neq 1$ we must have $c(f^*\xi)\neq 1$. \qed

If $X$ is a connected closed orientable manifold we have a degree one map $f:X\longrightarrow S^{\mathrm{dim}(X)}$. This induces an isomorphism in top integral cohomology. Thus we have the following. 

\begin{Cor}\label{secondcor}
Suppose $X$ is a connected closed orientable manifold. Then $\Sigma^k X$ is not $C$-trivial whenever $\mathrm{dim}(X)+k$ is even. \qed
\end{Cor}

If $k$ is odd then the suspensions $\Sigma^k (\mathbb C\mathbb P^m/\mathbb C\mathbb P^n)$  and $\Sigma^k (\mathbb H\mathbb P^m/\mathbb H\mathbb P^n)$ have cells only in odd dimensions and hence are $C$-trivial. The following observation is now immediate from the the above noted facts.  

\begin{Cor}\label{complex}
Let $\mathbb F= \mathbb C$ or $\mathbb H$. Let $0\leq n<m$. Then $\Sigma^k(\mathbb F\mathbb P^m/\mathbb F\mathbb P^n)$ is $C$-trivial if and only if $k$ is odd. \qed
\end{Cor}




\begin{Cor}
The product of two connected closed orientable manifolds is not $C$-trivial.
\end{Cor}
{\bf Proof.} In the case the product $M\times N$ is even dimensional, the claim follows from Corollary\,\ref{secondcor}. In the case that $M\times N$ is odd dimensional, assume that $M$ is even dimensional. Then $M$ is not $C$-trivial and since the composition  
$$M\stackrel{i}\longrightarrow M\times N\stackrel{\pi_1}{\longrightarrow} M$$
where $i(x)=(x,y)$ for a fixed $y\in N$ and $\pi_1$ the projection to the first factor, is identity it follows that $M\times N$ is not $C$-trivial. This completes the proof. \qed

In particular, a product of spheres is not $C$-trivial. We make some more observations before proving the main theorems. 

\begin{Lem}\label{howdoyoudo} Suppose $k$ is even. 
\begin{enumerate}
\item If $X^k_n$ is $W$-trivial, then $X^k_n$ is $C$-trivial. 
\item If $n$ is even and $X^k_{m,n}$ is $W$-trivial, then $X^k_{m,n}$ is $C$-trivial.
\end{enumerate}
\end{Lem} 
{\bf Proof.} We prove (1), the proof of (2) is similar. As $k$ is even, the integral cohomology is zero in odd degrees except in degree $n$ when $n$ is odd in which case it is infinite cyclic. 
The integral cohomology in even degrees is cyclic of order two. The mod-$2$ reduction map in even degrees is readily seen to be an isomorphism. 
By Lemma\,\ref{prelim} (4), the proof is complete. \qed


\begin{Lem}\label{newlemma}
Let $k,m,n$ be even. If $X^k_{m+1,n}$ is not $C$-trivial, then $X^k_{m,n}$ is not $C$-trivial. 
\end{Lem}
{\bf Proof.} The lemma follows from the fact that the obvious map 
$$j:X_{m,n}\longrightarrow X_{m+1,n}$$
induces isomorphism in integral cohomology in even degrees. Hence so does the map $\Sigma^kj$. \qed

{\em Proof of Theorem\,\ref{secondtheorem}.}  
If both $k$ and $n$ are odd, then as $\mathbb R\mathbb P^n$ is orientable, it follows from 
Corollary\,\ref{secondcor} that $X^k_n$ is not $C$-trivial.  Next assume that 
$k$ is odd and $n$ is even. Then as $H^i(X^k_n;\mathbb Z)=0$ if $i>0$ is even it follows that $X^k_n$ is $C$-trivial in this case. This proves the theorem when $k$ is odd.

Next we assume that $k$ is even. We first show that $X^4_n$ is not $C$-trivial if and only if  $n\geq 4$.
We know that $X^4_n$ for $n \leq 3$ is $W$-trivial (\cite{tanaka}, Theorem\,1.4). 
Hence it follows from Lemma\,\ref{howdoyoudo} that $X^4_n$ is $C$-trivial for $n \leq 3$. So assume that $n\geq 4$. Let $\xi$ be a complex $2$-plane bundle over $S^4$ with $c_2(\xi)$ a generator and let $\eta$ be a complex line bundle over $\mathbb R\mathbb P^n$ with 
$c_1(\eta)=t\in H^2(\mathbb R\mathbb P^n;\mathbb Z)\cong \mathbb Z_2$ the non-zero element. 
The cofiber sequence 
$$S^4\vee \mathbb R\mathbb P^n\stackrel{j}\longrightarrow S^4\times \mathbb R\mathbb P^n\stackrel{\alpha}\longrightarrow\Sigma^4\mathbb R\mathbb P^n$$
gives rise to an exact sequence 
$$0\rightarrow \widetilde{K}(\Sigma^4\mathbb R\mathbb P^n)\stackrel{\alpha}\longrightarrow\widetilde{K}(S^4\times\mathbb R\mathbb P^n)
\stackrel{j}\longrightarrow\widetilde{K}(S^4\vee \mathbb R\mathbb P^n)\rightarrow 0.$$
We compute (see, for example, Lemma\,2.1, \cite{tanaka})
$$\begin{array}{rcl}
c((\pi_1^*\xi-2)\otimes (\pi_2^*\eta-1)) & = & 1+ c_2(\xi)\times ((1+t)^{-2}-1)\\
& = & 1+ c_2(\xi)\times (-2t+3t^2-4t^3+\cdots)\\
& \neq & 1
\end{array}$$
as $c_2(\xi)\neq 0$ is a generator and $3t^2\neq 0$. 
Now as $j((\pi_1^*\xi-2)\otimes (\pi_2^*\eta-1))=0$, there exists a bundle $\theta\in \widetilde{K}(\Sigma^k\mathbb R\mathbb P^n)$ with 
$$\alpha(\theta)=(\pi_1^*\xi-2)\otimes (\pi_2^*\eta-1).$$ 
Clearly, $c(\theta)\neq 1$. This completes the proof that $X^4_n$ is not $C$-trivial if and only if $n\geq 4$. 

We now look at the case $k = 2$. Note that $X^2_1=S^3$ is $C$-trivial. We now check that $X^2_n$ is not $C$-trivial if $n \geq 2$. 
Let $\xi$ be a complex bundle over $S^{2}$ with $c(\xi)\neq 1$ and $c_1(\xi)$ a generator. Let $\eta$ denote the non-trivial complex line bundle over $\mathbb R\mathbb P^n$. Let $\pi_1,\pi_2$ be the two projections of $S^{2}\times\mathbb R\mathbb P^n$ onto the the first and the second factor respectively. Then, 
$$\begin{array}{rcl}
c((\pi_1^*\xi-1)\otimes (\pi_2^*\eta-1)) & = & 1+ c_1(\xi)\times ((1+t)^{-1}-1)\\
& = & 1+ c_1(\xi)\times (-t+t^2-t^3+\cdots)\\
& \neq & 1
\end{array}$$
as $c_1(\xi)$ is a generator. 
Here $t\in H^2(\mathbb R\mathbb P^n;\mathbb Z)=\mathbb Z_2$ is the unique non-zero element. Then, arguing as in the above case, it follows that 
there must exist a bundle $\theta\in \widetilde{K}(\Sigma^k\mathbb R\mathbb P^n)$ with $c(\theta)\neq 1$. Thus $X^2_n$ is not $C$-trivial if and only if  $n\geq 2$.





To complete the proof of the theorem we finally show that $X^6_n$ is $C$-trivial for all $n>0$. This will imply, by Theorem\,\ref{hahaha}, that $X^k_n$ is $C$-trivial for all $k\geq 6$ and $k$ even.  
First note that $X^6_1=S^7$ is $C$-trivial. That $X^6_2$ is $C$-trivial follows from the fact that $X^4_2$ is $C$-trivial and by Theorem\,  \ref{hahaha}.
By Theorem\,1.4 of \cite{tanaka}, $X^6_n$ is $W$-trivial whenever $n>3$. 
Hence by Lemma\,\ref{howdoyoudo} $X^6_n$ is $C$-trivial when $n>3$. Finally, we look at $X^6_3$. The long exact sequence of the pair $(\mathbb R\mathbb P^3,\mathbb R\mathbb P^2)$ shows that 
the inclusion map $i:\mathbb R\mathbb P^2\longrightarrow \mathbb R\mathbb P^3$ induces an isomorphism $i:H^2(\mathbb R\mathbb P^3;\mathbb Z)\longrightarrow H^2(\mathbb R\mathbb P^2;\mathbb Z)$. Hence the map $\Sigma^6i:\Sigma^6\mathbb R\mathbb P^2\longrightarrow \Sigma^6\mathbb R\mathbb P^3$ induces isomorphism in integral cohomology in degree $8$. Since the only non-zero cohomology in even degree (for both $X^6_2$ and $X^6_3$) is in 
degree $8$ and $X^6_2$ is $C$-trivial, it follows by Lemma\,\ref{prelim1} (2) that $X^6_3$ is $C$-trivial. This completes the proof that $X^6_n$ is $C$-trivial.    

This takes care of all the cases and completes the proof of the theorem. \qed

We now come to the proof of Theorem\,\ref{thirdtheorem}. First note that if $m$ is odd then the stunted real projective space $X_{m,m-2}$ 
admits a splitting 
$$X_{m,m-2}=S^m\vee S^{m-1}$$
and if $m$ is even then 
$$X_{m,m-2}=\Sigma^{m-2}\mathbb R\mathbb P^2.$$ 

We now prove Theorem\,\ref{thirdtheorem}.

{\em Proof of Theorem\,\ref{thirdtheorem}.} We first prove (1). If $k$ is odd and $m$ is even, then the integral cohomology of $X^k_{m,n}$ is trivial in even degrees and hence $X^k_{m,n}$ is $C$-trivial in this case proving (1). 

Next, if $k,m$ are both odd there exists a map $\alpha:X^k_{m,n}\longrightarrow S^{k+m}$ inducing isomorphism in top integral cohomology. By Corollary\,\ref{secondcor}, $X^k_{m,n}$ is not $C$-trivial. This proves (2). 

Next we prove (3). By Theorem 7.3 \cite{adams}, the projection $j:\mathbb R\mathbb P^m \rightarrow X_{m,n}$ maps $\widetilde{K}(X_{m,n})$ isomorphically into the subgroup of $\widetilde{K}(\mathbb R\mathbb P^m)$ generated by class of $2^t\nu$, where $\nu = \xi \otimes \mathbb C$ is the complexification of the canonical line bundle $\xi$ over $\mathbb R \mathbb P^m$. Let $\alpha \in \widetilde{K}(X_{m,n})$ be the generator such that $j^*(\alpha) = 2^t\nu$.
If $z  \in H^2(\mathbb R \mathbb P^m;\mathbb Z)$ is the unique non-zero element, then the total Chern class
$$c(2^t\nu) = c(\nu)^{2^t} = (1+z)^{2^t} = (1+z^{2^t}).$$
Since $j^*:H^k(X_{m,n};\mathbb Z) \rightarrow H^k(\mathbb R \mathbb P^m;\mathbb Z)$ is injective  for $0 \leq k \leq m$, we have $c(\alpha) = 1$ if and only if $m < 2^{t+1}$. As $\alpha$ is a generator, we conclude that $X_{m,n}$ is $C$-trivial if and only of $m< 2^{t+1}$.




We now prove (4).
We shall only prove the $C$-triviality for the case $k=2$. Then in view of Theorem\,\ref{hahaha}, $X^k_{m,n}$ will be  $C$-trivial for all $k\geq 4$ and $k$ is even.  By Theorem\, 1.3   of \cite{aniajay}, $X^2_{m,n}$ is $W$-trivial if $m \neq 6,7$. By Lemma\,\ref{howdoyoudo} we see that for $n$ even, $X^2_{m,n}$ is $C$-trivial if $m \neq 6,7$. We are now left to prove $C$-triviality of the following cases: $X^2_{6,n}$ and $X^2_{7,n}$ for $n$ even. We prove these as follows.

We first prove that $X^2_{6,n}$ is $C$-trivial for $n$ even. As $X^2_{6,4}=\Sigma^6\mathbb R\mathbb P^2$, it follows from Theorem\,\ref{secondtheorem} that 
$X^2_{6,4}$ is $C$-trivial. We next look at the case $X^2_{6,2}$ and  
let $\xi$ be a complex vector bundle over $X^2_{6,2}$. We shall show that $c_3(\xi)=0=c_4(\xi)$. 
We first claim that $c_3(\xi)=0$. For if $c_3(\xi)\neq 0$, then $w_6(\xi_{\mathbb R})\neq 0$ as the mod-$2$ reduction homomorphism is an isomorphism. 
Now observe that $w_i(\xi_{\mathbb R})=0$ for $1\leq i\leq 5$. This contradicts the well-known fact that for a real bundle the first non-zero Stiefel-Whitney class appears in degree a power of two. Thus $c_3(\xi)=0$. Assume now that $c_4(\xi)\neq 0$. 
consider the exact sequence 
$$\cdots\rightarrow \widetilde{K}^{-2}(X_{8,6})\stackrel{\alpha}\longrightarrow\widetilde{K}^{-2}(X_{8,2})\stackrel{j}\longrightarrow
\widetilde{K}^{-2}(X_{6,2})\longrightarrow\widetilde{K}^{-1}(X_{8,6})\rightarrow \cdots.$$
Using the Atiyah-Hirzebruch spectral sequence in complex $K$-theory it is easy to see that if  $m$ is even, then 
$$\widetilde{K}^{-1}(X_{m,n})=0.$$ 
Thus the homomorphism $j$ is epimorphic. Hence, $X^2_{6,2}$ is $C$-trivial as $X^2_{8,2}$ is $C$-trivial. 

Next we prove that $X^2_{7,n}$ is $C$-trivial for $n$ even.  Clearly, $X^2_{7,6}=S^9$ is $C$-trivial. 
That $X^2_{7,2}$ and $X^2_{7,4}$ is $C$-trivial follows from $C$-triviaility of $X^2_{6,2}$ and $X^2_{6,4}$ and by Lemma\,\ref{newlemma}. This completes the proof of (4).

We now prove (5). Here $k$ is even and $n$ is odd. First we assume $m$ is even. We look at the cofiber sequence 
$$X_{n+1,n}\stackrel{j}\longrightarrow X_{m,n}\stackrel{\alpha}\longrightarrow X_{m,n+1}$$
and the associated exact sequence
$$\cdots\rightarrow \widetilde{K}^{-k}(X_{m,n+1})\stackrel{\alpha}\longrightarrow\widetilde{K}^{-k}(X_{m,n})\stackrel{j}\longrightarrow
\widetilde{K}^{-k}(X_{n+1,n})\longrightarrow\widetilde{K}^{-k+1}(X_{m,n+1})\rightarrow\cdots.$$
As noted above, since the last group in the above exact sequence is zero the homomorphism $j$ is epimorphic. Since $\Sigma^kX_{n+1,n}$ is an even dimensional sphere, and therefore not $C$-trivial, we conclude that $X^k_{m,n}$ is not $C$-trivial. Note that, in particular, there exists a complex vector bundle $\xi$ over $X^k_{m,n}$ with $$c_{\frac{k+n+1}{2}}(\xi)\neq 0.$$ 

Next we assume $m$ is odd. Consider the obvious map 
$$j: X_{m,n}\longrightarrow X_{m+1,n}.$$
The homomorphism $j$ in integral cohomology is an isomorphism in degree $(n+1)$. The homomorphism $\Sigma^kj$ is an isomorphism in integral cohomology in degree $(k+n+1)$. Now if $\xi$ is a complex vector bundle over $X^k_{m+1,n}$ with $$c_{\frac{k+n+1}{2}}(\xi)\neq 0,$$ we have that $$c_{\frac{k+n+1}{2}}(j^*\xi)\neq 0.$$ 
Thus $X^k_{m,n}$ is not $C$-trivial when $k$ is even and $m,n$ are odd.  This completes the proof of (5) and the theorem. 
\qed

We remark that the fact that 
$$\widetilde{K}^{-k}(X_{m,n})=0$$
whenever $k$ is odd and $m$ is even also gives another proof of Theorem\,\ref{thirdtheorem} (1).


\begin{thebibliography}{99}

\bibitem{atiyahirz} Atiyah, M., and Hirzebruch, F., {\em Bott periodicity and the parallelizability of the spheres},  Proc. Cambridge Philos. Soc., 57 (1961)  223-226. 
\bibitem{adams} Adams, J. F., {\em Vector fields on spheres}, Ann. of Math., Second Series, Vol. 75, No. 3 (1962) 603-632. 


\bibitem{hus} Husemoller D., {\em Fibre Bundles}, Springer-Verlag, New York 1966. 
\bibitem{aniajay} Naolekar, A. C., and Thakur, A. S., {\em Vector bundles over iterated suspensions of stunted real projective spaces.}, Acta Math. Hungarica, to appear. Available at http://arxiv.org/pdf/1211.5274.
\bibitem{tanaka1} Tanaka, R., {\em On trivialities of Stiefel-Whitney classes of vector bundles over highly connected complexes}, Top. and App., 155 (2008)  1687-1693.
\bibitem{tanaka} Tanaka, R., {\em On trivialities of Stiefel-Whitney classes of vector bundles over iterated suspension spaces}, HHA, vol (12) No. 1 (2010) 357-366. 
\bibitem{ajay} Thakur, A. S., {\em On trivialities of Stiefel-Whitney classes of vector bundles over iterated suspensions of Dold manifolds}, HHA, 15 (2013), 223-233.



\end{thebibliography}
\end{document}